\documentclass{amsart} 
\usepackage{latexsym} 
\usepackage{amssymb} 
\usepackage{amsmath} 

\begin{document}


\newtheorem{theorem}{Theorem} 
\newtheorem{problem}{Problem} 
\newtheorem{definition}{Definition} 
\newtheorem{lemma}{Lemma} 
\newtheorem{proposition}{Proposition} 
\newtheorem{corollary}{Corollary} 
\newtheorem{example}{Example} 
\newtheorem{conjecture}{Conjecture} 
\newtheorem{algorithm}{Algorithm} 
\newtheorem{exercise}{Exercise} 
\newtheorem{remarkk}{Remark} 
 
\newcommand{\be}{\begin{equation}} 
\newcommand{\ee}{\end{equation}} 
\newcommand{\bea}{\begin{eqnarray}} 
\newcommand{\eea}{\end{eqnarray}} 
\newcommand{\beq}[1]{\begin{equation}\label{#1}} 
\newcommand{\eeq}{\end{equation}} 
\newcommand{\beqn}[1]{\begin{eqnarray}\label{#1}} 
\newcommand{\eeqn}{\end{eqnarray}} 
\newcommand{\beaa}{\begin{eqnarray*}} 
\newcommand{\eeaa}{\end{eqnarray*}} 
\newcommand{\req}[1]{(\ref{#1})} 
 
\newcommand{\lip}{\langle} 
\newcommand{\rip}{\rangle} 

\newcommand{\uu}{\underline} 
\newcommand{\oo}{\overline} 
\newcommand{\La}{\Lambda} 
\newcommand{\la}{\lambda} 
\newcommand{\eps}{\varepsilon} 
\newcommand{\om}{\omega} 
\newcommand{\Om}{\Omega} 
\newcommand{\ga}{\gamma} 
\newcommand{\rrr}{{\Bigr)}} 
\newcommand{\qqq}{{\Bigl\|}} 
 
\newcommand{\dint}{\displaystyle\int} 
\newcommand{\dsum}{\displaystyle\sum} 
\newcommand{\dfr}{\displaystyle\frac} 
\newcommand{\bige}{\mbox{\Large\it e}} 
\newcommand{\integers}{{\Bbb Z}} 
\newcommand{\rationals}{{\Bbb Q}} 
\newcommand{\reals}{{\rm I\!R}} 
\newcommand{\realsd}{\reals^d} 
\newcommand{\realsn}{\reals^n} 
\newcommand{\NN}{{\rm I\!N}} 
\newcommand{\DD}{{\rm I\!D}} 
\newcommand{\degree}{{\scriptscriptstyle \circ }} 
\newcommand{\dfn}{\stackrel{\triangle}{=}} 
\def\complex{\mathop{\raise .45ex\hbox{${\bf\scriptstyle{|}}$} 
     \kern -0.40em {\rm \textstyle{C}}}\nolimits} 
\def\hilbert{\mathop{\raise .21ex\hbox{$\bigcirc$}}\kern -1.005em {\rm\textstyle{H}}} 
\newcommand{\RAISE}{{\:\raisebox{.6ex}{$\scriptstyle{>}$}\raisebox{-.3ex} 
           {$\scriptstyle{\!\!\!\!\!<}\:$}}} 
 
\newcommand{\hh}{{\:\raisebox{1.8ex}{$\scriptstyle{\degree}$}\raisebox{.0ex} 
           {$\textstyle{\!\!\!\! H}$}}} 

\newcommand{\OO}{\won} 
\newcommand{\calA}{{\mathcal A}} 
\newcommand{\calB}{{\mathcal B}} 
\newcommand{\calC}{{\cal C}} 
\newcommand{\calD}{{\cal D}} 
\newcommand{\calE}{{\cal E}} 
\newcommand{\calF}{{\mathcal F}} 
\newcommand{\calG}{{\cal G}} 
\newcommand{\calH}{{\cal H}} 
\newcommand{\calK}{{\cal K}} 
\newcommand{\calL}{{\mathcal L}} 
\newcommand{\calM}{{\cal M}} 
\newcommand{\calO}{{\cal O}} 
\newcommand{\calP}{{\cal P}} 
\newcommand{\calU}{{\mathcal U}} 
\newcommand{\calX}{{\cal X}} 
\newcommand{\calXX}{{\cal X\mbox{\raisebox{.3ex}{$\!\!\!\!\!-$}}}} 
\newcommand{\calXXX}{{\cal X\!\!\!\!\!-}} 
\newcommand{\gi}{{\raisebox{.0ex}{$\scriptscriptstyle{\cal X}$} 
\raisebox{.1ex} {$\scriptstyle{\!\!\!\!-}\:$}}} 
\newcommand{\intsim}{\int_0^1\!\!\!\!\!\!\!\!\!\sim} 
\newcommand{\intsimt}{\int_0^t\!\!\!\!\!\!\!\!\!\sim} 
\newcommand{\pp}{{\partial}} 
\newcommand{\al}{{\alpha}} 
\newcommand{\sB}{{\cal B}} 
\newcommand{\sL}{{\cal L}} 
\newcommand{\sF}{{\cal F}} 
\newcommand{\sE}{{\cal E}} 
\newcommand{\sX}{{\cal X}} 
\newcommand{\R}{{\rm I\!R}} 
\renewcommand{\L}{{\rm I\!L}} 
\newcommand{\vp}{\varphi} 
\newcommand{\N}{{\rm I\!N}} 
\def\ooo{\lip} 
\def\ccc{\rip} 
\newcommand{\ot}{\hat\otimes} 
\newcommand{\rP}{{\Bbb P}} 
\newcommand{\bfcdot}{{\mbox{\boldmath$\cdot$}}} 
 
\renewcommand{\varrho}{{\ell}} 
\newcommand{\dett}{{\textstyle{\det_2}}} 
\newcommand{\sign}{{\mbox{\rm sign}}} 
\newcommand{\TE}{{\rm TE}} 
\newcommand{\TA}{{\rm TA}} 
\newcommand{\E}{{\rm E\,}} 
\newcommand{\won}{{\mbox{\bf 1}}} 
\newcommand{\Lebn}{{\rm Leb}_n} 
\newcommand{\Prob}{{\rm Prob\,}} 
\newcommand{\sinc}{{\rm sinc\,}} 
\newcommand{\ctg}{{\rm ctg\,}} 
\newcommand{\loc}{{\rm loc}} 
\newcommand{\trace}{{\,\,\rm trace\,\,}} 
\newcommand{\Dom}{{\rm Dom}} 
\newcommand{\ifff}{\mbox{\ if and only if\ }} 
\newcommand{\nproof}{\noindent {\bf Proof:\ }} 
\newcommand{\remark}{\noindent {\bf Remark:\ }} 
\newcommand{\remarks}{\noindent {\bf Remarks:\ }} 
\newcommand{\note}{\noindent {\bf Note:\ }}

\newcommand{\boldx}{{\bf x}} 
\newcommand{\boldX}{{\bf X}} 
\newcommand{\boldy}{{\bf y}} 
\newcommand{\boldR}{{\bf R}} 
\newcommand{\uux}{\uu{x}} 
\newcommand{\uuY}{\uu{Y}} 
 
\newcommand{\limn}{\lim_{n \rightarrow \infty}} 
\newcommand{\limN}{\lim_{N \rightarrow \infty}} 
\newcommand{\limr}{\lim_{r \rightarrow \infty}} 
\newcommand{\limd}{\lim_{\delta \rightarrow \infty}} 
\newcommand{\limM}{\lim_{M \rightarrow \infty}} 
\newcommand{\limsupn}{\limsup_{n \rightarrow \infty}} 
 
\newcommand{\ra}{ \rightarrow }

\newcommand{\ARROW}[1] 
  {\begin{array}[t]{c}  \longrightarrow \\[-0.2cm] \textstyle{#1} \end{array} } 
 
\newcommand{\AR} 
 {\begin{array}[t]{c} 
  \longrightarrow \\[-0.3cm] 
  \scriptstyle {n\rightarrow \infty} 
  \end{array}} 
 
\newcommand{\pile}[2] 
  {\left( \begin{array}{c}  {#1}\\[-0.2cm] {#2} \end{array} \right) } 
 
\newcommand{\floor}[1]{\left\lfloor #1 \right\rfloor} 
 
\newcommand{\mmbox}[1]{\mbox{\scriptsize{#1}}} 
 
\newcommand{\ffrac}[2] 
  {\left( \frac{#1}{#2} \right)} 
 
\newcommand{\one}{\frac{1}{n}\:} 
\newcommand{\half}{\frac{1}{2}\:} 
 
\def\le{\leq} 
\def\ge{\geq} 
\def\lt{<} 
\def\gt{>} 
 
\def\squarebox#1{\hbox to #1{\hfill\vbox to #1{\vfill}}} 
\newcommand{\nqed}{\hspace*{\fill} 
           \vbox{\hrule\hbox{\vrule\squarebox{.667em}\vrule}\hrule}\bigskip} 
 
\title{Transportation cost inequalities for diffusions under uniform distance}

\author{ A. S. \"Ust\"unel} 
\maketitle 
\noindent 
{\bf Abstract:}{\small{ We prove the transportation inequality with
    the uniform norm  for the
    laws of diffusion processes with Lipschitz and/or dissipative
    coefficients and apply them to some singular stochastic
    differential equations of interest.
}}\\ 

\vspace{0.5cm} 

\section{\bf{Introduction} }
\noindent
Let $(W,d)$ be a separable Fr\'echet  space, for two probability measures
$P$and $Q$ on $(W,\calB(W))$, then  the Wasserstein distance
(cf. \cite{Vil})  between
$P$ and $Q$, denoted as $d_W(P,Q)$, is defined as 
$$
d_W^2(P,Q)=\inf\left\{\int_{W\times
    W}d(x,y)^2\theta(dx,dy):\,\theta\in \Sigma(P,Q)\right\}\,,
$$
where $\Sigma(P,Q)$ denotes the set of probability measures on
$W\times W$ whose first marginal is $P$ and the second one is $Q$;
note that this is a compact set under the weak topology, hence the
infimum is always attained for any $d$  (even lower
semi-continuous). It is quite useful to find an upper bound for this
distance, if possible  dimension independent. There are a lot of works
on this subject (cf. \cite{Vil}), beginning by the contributions of M. Talagrand,
cf. \cite{Tal}, where it is shown that the relative entropy is a  fully
satisfactory  upper bound. In \cite{fandu,fandu1}, it is shown that the
relative entropy is again an upper bound when $P$ is the Wiener measure
and $d$ is the singular Cameron-Martin distance using the Girsanov
theorem (cf. also \cite{DJ}). The same method has also been employed
in \cite{wz1}  and more recently  in
\cite{SPAL}  to obtain 
a transportation cost inequality w.r. to Banach norm for  diffusion
processes. The former assumes quite strong conditions on the
coefficients which govern the diffusion which are superfluous and make
difficult the applicability of the inequality, while the latter one
treats essentially the 
one-dimensional case with an extension 
 to the case where the diffusion coefficients
are independent and their slight perturbations. Inspired with these 
works, we have attacked the general case: namely, the case of fully
dependent diffusion like processes and their extensions and infinite
dimensional diffusion processes governed with a cylindrical Brownian
motion. Besides, there is a special class of diffusion processes  with
singular (dissipative)  drifts which are constructed as weak limits of
the Lipschitzian case where the approximating diffusions have
Lipschitz continuous drifts but the Lipschitz constant explodes at the
limit; this last class  is particularly interesting because of their
applications to physics. 

\noindent
To achieve this program, we need the following result about the
stability of the transportation cost inequality under the weak limits 
of probability measures,  which is proved by Djellout, Guillin and  Wu
in \cite{DJ}.  Since we make an important use of it, we  give it  with
a (slightly different and more general) proof.

\begin{lemma}
\label{stab-lemma}
Assume that $(P_k,k\geq 1)$ is a sequence of probability measures on a
separable Fr\'echet space $(W,d)$, converging weakly to a probability
$P$. If 
$$
d_W^2(Q,P_k)\leq c_k \int_W \frac{dQ}{dP_k}\log \frac{dQ}{dP_k} dP_k=c_kH(Q|P_k)
$$
for any $k\geq 1$, for any probability $Q$,  where $c_k>0$ are bounded
constants,  then the transportation 
inequality holds for $P$, namely
\begin{equation}
\label{tal-ineq}
d_W^2(Q,P)\leq c H(Q|P)\,,
\end{equation}
where $c=\sup_kc_k$.
\end{lemma}
\nproof
 If $f=dQ/dP$ is a bounded, continuous function, then the inequality
(\ref{tal-ineq}) follows from the lower semi continuity of the
transportation cost w.r. to the weak convergence  and from the
hypothesis since $f\log f$ is continuous and bounded. Due to the
dominated convergence theorem, to prove the general case, it suffices
to prove the case where $f$ is $P$-essentially bounded and
measurable. In this case, there exists a sequence of bounded, upper
semi continuous functions, say $(f_n,n\geq 1)$, increasing to $f$
$P$-almost surely. By the dominated convergence theorem, the measures 
$(\tilde{f}_ndP,\,n\geq 1)$ converge weakly  to the measure $fdP$,
where $\tilde{f}_n=f/P(f_n)$. On the other hand
$H(\tilde{f_n}dP|dP)\to H(fdP|P)$ again  by the dominated convergence
theorem. Hence, to prove the general case, it is sufficient to prove
the inequality with $f$ upper semi continuous and bounded. Since we
are on a Fr\'echet space, there exists a sequence of (positive)
continuous functions decreasing to $f$ which may be chosen uniformly
bounded by taking the minimum of each with the upper bound of f,  and the inequality
(\ref{tal-ineq}) follows again due to the dominated convergence
theorem.
\nqed

\section{\bf{Diffusion type processes with Lipschitz coefficients}}
\noindent
Let $(W,H,\mu)$ be the classical Wiener space, i.e.,
$W=C_0([0,1],\R^d),\,H=H^1([0,1],\R^d)$ and $\mu$ is the Wiener
measure under which the evaluation map at $t\in [0,1]$ is a Brownian
motion. Suppose that $X=(X_t,t\in[0,1])$ is the solution of the
following SDE (stochastic differential equation)
\beaa
dX_t&=&\sigma(t,X_t)dW_t+b(t,X)dt\\
X_0&=&z\in\R^d
\eeaa
where $\sigma:[0,1]\times \R^d\to\otimes\R^d$ is  uniformly
Lipschitz w.r.to $x$ with a Lipschitz constant being equal to $K$, $b:[0,1]\times W\to \R^d$ is adapted and such that 
$$
|b(t,\xi)-b(t,\eta)|\leq K\sup_{s\leq
  t}|\xi(s)-\eta(s)|=\|\xi-\eta\|_{t}
$$
for any $\xi,\eta\in W$. We denote by $d_W$ the Wasserstein distance
on the probability measures on $W$ defined by the uniform norm:
$$
d^2_W(\rho,\nu)=\inf\left(\int_{W\times
    W}\|x-y\|^2d\ga(x,y):\,\ga\in\Sigma(\rho,\nu)\right)
$$
where $\Sigma(\rho,\nu)$ the set of probabilities on $W\times W$ whose
first marginals are $\rho$ and the second one is $\nu$. We have the
following bound for $d_W$:
\begin{theorem}
\label{th-1}
Let $P$ be the law of the solution of the SDE described above, then
for any probability $Q$ on $(W,\calB(W))$, we have
\begin{equation}
\label{ineq-1}
d^2_W(P,Q)\leq 6\,e^{15K^2}H(Q|P)\,
\end{equation}
where $H(Q|P)$ is the relative entropy of $Q$ w.r. to $P$.
\end{theorem}
\nproof
Due to the rotation invariance of the Wiener measure, we can suppose without loss of generality that $\sigma$ takes its
values in the set of positive matrices. Suppose first that $\sigma$ is
strictly elliptic. From the general results about the SDE (cf. \cite{I-W,R-W}), the
coordinate process $x$ under
the probability $P$ can be written as 
$$
dx_t=\sigma(t,x_t)d\beta_t+b(t,x)dt
$$
with $x_0=z$ $P$-a.s., where $\beta$ is an $\R^d$-valued $P$-Brownian
motion. At this point of the proof we need the following result, which
is probably well-known (cf. \cite{R-W} and the references there),
though we include its  proof for the sake of completeness:
\begin{lemma}
\label{lem-1}
Any bounded  $P$-martingale can be written as a stochastic integral w.r. to $\beta$ of an
adapted process $(\alpha_s,s\in [0,1])$, with $E_P\int_0^1|\alpha_s|^2ds<\infty$.
\end{lemma}
\nproof
Let us denote by $P^0$ the law of the solution of 
$$
dX_t=\sigma (t,X_t)dW_t\,,
$$
then under $P^0$, the coordinate process $x$ can be written as 
$$
dx=\sigma(t,x_t)d\beta^0_t\,,
$$
where $\beta^0$ is a $P^0$-Brownian motion.
Let $Z$ be a bounded $P$-martingale with $Z_0=0$, assume that it is orthogonal
to the Hilbert  space of $P$-square integrable martingales written as the stochastic integrals
w.r. to $\beta$ of the adapted processes. Let $M$ be the exponential martingale
defined as 
$$
M_t=\exp\left(-\int_0^t(\sigma^{-1}(s,x_s)b(s,x),d\beta_s)-\frac{1}{2}\int_0^t|\sigma^{-1}(s,x_s)b(s,x)|^2ds\right)\,.
$$
Then, we know from the uniqueness and the Girsanov theorem  that
$MdP=dP^0$, since $M$ can be written as a stochastic integral w.r. to
$\beta$, our hypothesis implies that $ZM$ is again a $P$-martingale,
hence $Z$ is a $P^0$-martingale, therefore, from the classical Markov
case it can be written as 
\beaa
Z_t&=&\int_0^tH_s.d\beta^0_s\\
&=&\int_0^tH_s.(d\beta_s-\sigma^{-1}(s,x_s)b(s,x)ds)\,.
\eeaa
This last expression implies that 
$$
\langle Z,Z\rangle_t=\langle Z,\int_0^\cdot H_s.d\beta_s\rangle_t\,
$$
but $Z$ is orthogonal to the stochastic integrals of the form $\int
\alpha_s.d\beta_s$, hence $Z_t=E_P[Z_t]=0$, which proves the claim.
\nqed

\noindent
Let us complete now the proof of the theorem: If $Q$ is singular w.r. to $P$, then there is nothing to prove
due to the definition of the entropy. Let $L$ be the Radon-Nikodym
derivative $dQ/dP$, we shall first  suppose that $L>0$ $P$-a.s. In
this case we can write
$$
L=\rho(-\delta v)\,,
$$
where $v(t,x)=\int_0^t\dot{v}_s(x)ds$, $\dot{v}_s(x)$ is a.s. adapted
and $\int_0^1|\dot{v}_s(x)|^2ds<\infty$ a.s. and $\delta
v=\int_0^1\dot{v}_sd\beta_s$. From the Girsanov theorem,
$z_t=\beta_t+\int_0^t\dot{v}_sds$ is $Q$-Brownian motion, hence by the
uniqueness of the solution of SDE, if we denote by $x^v$ the solution
of the SDE given as 
$$
dx^v_t=\sigma(t,x^v_t)dz_t+b_t(x^v)dt\,
$$
the image of $Q$ under the solution map $x^v$ is equal to $P$,
consequently $(x^v\times I_W)(Q)\in \Sigma(P,Q)$, hence we have the
following domination:
$$
d_W^2(P,Q)\leq E_Q[\|x^v-x\|^2]
$$
where $\|\cdot\|$ denotes the uniform norm on $W$. Using Doob and
H\"older inequalities, we get
\beaa
E_Q[\sup_{r\leq t}|x^v_r-x_r|^2] &\leq&(12+3t)K^2E_Q\int_0^t|x^v_s-x_s|^2ds\\
&&+3t E_Q\int_0^t|\dot{v}_s|^2ds\,.
\eeaa
It follows from the Gronwall lemma that 
$$
E_Q[\sup_{r\leq t}|x^v_r-x_r|^2]\leq
3t\,E_Q\int_0^t|\dot{v}_s|^2ds\,e^{3K^2(4+t)}\,
$$
since 
$$
E_Q\int_0^1|\dot{v}_s|^2ds=2H(Q|P)
$$
the claim follows in the case $P\sim Q$. For the case where $Q\ll P$
let 
$$
L_\eps=\frac{L+\eps}{1+\eps}\,,
$$
then it is easy to see that $(L_\eps\log L_\eps,\,\eps\leq \eps_0)$ is
$P$-uniformly integrable provided $E_P[L\log L]<\infty$. Hence the
proof, in the strictly elliptic case, 
follows by the lower semi-continuity of $Q\to d_W(P,Q)$. The general
case follows by replacing $\sigma$ by $\eps I_{\R^d}+\sigma$, then
remarking that the corresponding probabilities $(P_\eps, \eps\leq
\eps_0)$ converge weakly and that 
$$
d_W^2(P_\eps,Q)\leq 6\,e^{15(\eps+K)^2}H(Q|P_\eps)\,
$$
and hence it follows from Lemma \ref{stab-lemma} that 
$$
d_W^2(P,Q)\leq 6\,e^{15K^2}H(Q|P)\,.
$$
\nqed

\noindent 
Since the inequality (\ref{ineq-1}) is dimension independent, we can
extend it easily to the infinite dimensional case:
\begin{corollary}
\label{infdim-thm}
Let $M$ be a separable Hilbert space, suppose that $B$ is a
$M$-cylindrical Wiener process. Assume that $\sigma: [0,1]\times M\to
L_2(M,KM=M\otimes_2M$ (space of Hilbert-Schmidt operators on $M$) and
$b:[0,1]\times M\to M$ are  uniformly Lipschitz with Lipschitz
constant $K$. Let $P$ be the law of the following SDE:
$$
dX_t=\sigma (t,X_t)dB_t+b(t,X_t)dt\,,X_0=x\in M\,.
$$
Then the law of $P$ satisfies the transportation cost inequality
(\ref{ineq-1}).
\end{corollary}
\nproof
Let $(\pi_n,n\geq 1)$ be an sequence of orthogonal projections of $M$
increasing to the identity, define $\sigma_n=\pi_n\sigma\circ \pi_n$,
$b_n=\pi_n b\circ \pi_n$, $B^n=\pi_nB$ and $x^n=\pi_n x$. Let then
$P^n $ be the law of the SDE 
$$
dX^n_t=\sigma^n (t,X^n_t)dB^n_t+b^n(t,X^n_t)dt\,,X^n_0=x^n\,.
$$
From Theorem \ref{th-1}, $P^n$ satisfies the inequality
(\ref{ineq-1}) with a  constant independent of $n$, since 
$(P^n,n\geq 1)$ converges weakly to $P$, the proof follows from Lemma \ref{stab-lemma}.
\nqed

\section{\bf{Transport inequality for the monotone case}}

\noindent
Assume  that the Lipschitz property of the adapted  drift coefficient is
replaced by the following dissipativity hypothesis 
$$
(b(t,x)-b(t,y),x_t-y_t)\leq 0
$$
for any $t\in[0,1]$ and $x,y\in W$, where, as before $(\cdot,\cdot)$
denotes the scalar product in $\R^d$. The derivative of a proper
concave function on $\R^d$ is a typical example of such
drift. We shall suppose first that 
$$
\int_0^1|b(t,x)|^2ds<\infty
$$
for any $x\in W$. 
\begin{proposition}
\label{prop-1}
Assume that $b$ is of linear growth, i.e., $|b(t,x)|\leq N(1+\|x\|)$
and let  $P$ be  the law of the solution of the following SDE
\begin{equation}
\label{msde}
dX_t=\sigma(t,X_t)dW_t+b(t,X)dt+m(t,X_t)dt
\end{equation}
with $X_0=x\in \R^d$ and that $\sigma$ and $m:[0,1]\times \R^d\to
\R^d$ are  uniformly $K$-Lipschitz w.r. to
the space variable.  Then for
any $Q\ll P$, we have 
\bea
\label{tr-2}
d_W^2(P,Q)&\leq &\left(c2^{3/2} \|\sigma\|_\infty^{3/2}
e^{\frac{1}{2}(K^2+2K+1)}\right)\sqrt{H(Q|P)}\\
&&+2\|\sigma\|_\infty e^{\frac{1}{2}(K^2+2K+1)}\left(1+K(K+2))e^{\frac{1}{2}(K^2+2K+1)}\right)H(Q|P)\,,\nonumber
\eea
where $\|\sigma\|_\infty$ is a uniform bound for $\sigma$, $K$ is
the Lipschitz constant and $c$ is the universal constant of Davis'
inequality for $p=1$.

\end{proposition}
\nproof 
Recall that under $P$, the coordinate process satisfies
$dx=\sigma(t,x_t)d\beta+(b(t,x)+m(t,x_t))dt$, where $\beta$ is a $P$-Brownian
motion. 
Assume that $Q$ is another probability on $W$ such that $Q\ll P$, let
$L$ be $dQ/dP$. Suppose first that $L>0$ $P$-almost surely. As
explained  in the first section, we can write $L$ as an exponential
martingale $L=\rho(-\delta v)$, then $x^v(Q)=P$, where $x^v$ is
defined as before: $dx^v=\sigma
(t,x^v_t)(d\beta_t+\dot{v}_tdt)+b(t,x^v)dt+m(t,x^v_t)dt$. Again by the uniqueness
of the solutions, we have $(x^v\times I_W)(Q)\in \Sigma(P,Q)$, hence 
$$
d_W^2(P,Q)\leq E_Q[\|x^v-x\|^2]\,.
$$
It follows from the It\^o formula, letting $dz=d\beta+\dot{v}dt$, that 
\beaa
|x^v_t-x_t|^2&=&2\int_0^t(x^v_s-x_s,dx_s^v-dx_s)+\int_0^t|\sigma(s,x^v_s)-\sigma(s,x_s)|^2ds\\
&=&2\int_0^t(x^v_s-x_s,b(s,x^v)-b(s,x))ds\\
&&+2\int_0^t(x^v_s-x_s,(\sigma(s,x^v_s)-\sigma(s,x_s))dz_s+(m(s,x^v_s)-m(s,x_s))ds)\\
&&+\int_0^t|\sigma(s,x^v_s)-\sigma(s,x_s)|^2ds-2\int_0^t(x^v_s-x_s,\sigma(s,x^v_s)\dot{v}_s)ds\,.
\eeaa
By the dissipative character of $b$, we get
\beaa
|x^v_t-x_t|^2&\leq&2\int_0^t(x^v_s-x_s,(\sigma(s,x^v_s)-\sigma(s,x_s))dz_s+(m(s,x^v_s)-m(s,x_s))ds)\\
&&+\int_0^t|\sigma(s,x^v_s)-\sigma(s,x_s)|^2ds-2\int_0^t(x^v_s-x_s,\sigma(s,x^v_s)\dot{v}_s))ds\,.
\eeaa
Using, the usual stopping techniques, we can suppose that the
stochastic integral has zero expectation and taking the
$Q$-expectation of both sides, we obtain
\beaa
E_Q[|x^v_t-x_t|^2]&\leq&(2K+K^2)E\int_0^t|x^v_s-x_s|^2ds\\
&&+2\|\sigma\|_\infty E\int_0^t|x^v_s-x_s||\dot{v}_s|ds
\eeaa
using the inequality $xy\leq \delta (x^2/2)+(y^2/2\delta)$, we get
$$
E_Q[|x^v_t-x_t|^2]\leq(2K+K^2+\delta \|\sigma\|_\infty^2)E\int_0^t|x^v_s-x_s|^2ds++\frac{2}{\delta}H_t(Q|P)\,,
$$
where $\delta>0$ is arbitrary and  $H_t(Q|P)=\int\log \frac{dQ}{dP}|_{\calF_t}dQ$ is  the entropy
for the horizon $[0,t]$, which is an increasing function of $t$. It
follows from the Gronwall lemma that 
\begin{equation}
\label{ineq-2}
E_Q[|x^v_t-x_t|^2]\leq
\frac{2}{\delta}H_t(Q|P)\exp\left[ t(2K+K^2+\delta \|\sigma\|_\infty^2)\right]\,.
\end{equation}
Using now the Davis' inequality, the  Lipschitz property and the
boundedness of $\sigma$,  we get
\beaa
E[\sup_{r\leq t}|x^v_r-x_r|^2]&\leq&(2c\|\sigma\|_\infty+\sqrt{2}H_t(Q|P)^{1/2})E\left[\int_0^t|x^v_s-x_s|^2ds\right]^{1/2}\\
&&+K(K+2)E\int_0^t|x^v_s-x_s|^2ds\,,
\eeaa
where $c$ is the universal constant of Davis' inequality. Note that
the right hand side of the inequality (\ref{ineq-2}) is monotone
increasing in $t$, we insert it to the above inequality and minimize
it w.r. to $\delta$ for $t=1$  and the proof is completed.
\nqed

\noindent
In fact we have another version of the inequality (\ref{ineq-2}) in
the case where $\sigma$ is not bounded but still $K$-Lipschitz:
\begin{proposition}
\label{prop-2}
Assume that all the hypothesis of Proposition \ref{prop-1} are
satisfied except the boundedness of $\sigma$ which appears in the SDE (\ref{msde}), then we have the
following transportation cost inequality:
\begin{equation}
\label{ineq-3}
d_W^2(P,Q)\leq
H(Q|P)\frac{2}{(1-acK)^2}\exp\left(\frac{1}{1-acK}\left(\frac{cK}{a}+1-acK+2K+K^2\right)\right)
\end{equation}
where $P$ is the law of the SDE (\ref{msde}),  $Q$ is any other
probability and $a>0$ is arbitrary provided that $acK<1$.
\end{proposition}
\nproof
The proof is somewhat similar to the proof of Proposition
\ref{prop-1}: in fact we control uniformly  the stochastic integral
term in the It\^o development of $|x^v_t-x_t|^2$ as follows:
\beaa
\lefteqn{E\left[\sup_{r\leq
  t}\left|\int_0^r(x^v_s-x_s,(\sigma(s,x^v_s)-\sigma(s,x_s)dz_s)\right|\right]}\\
&\leq& c E\left[\left(\int_0^t|x^v_s-x_s|^2|\sigma(s,x^v_s)-\sigma(s,x_s)|^2ds\right)^{1/2}\right]\\
&\leq& cK E\left[\left(\int_0^t|x^v_s-x_s|^4ds\right)^{1/2}\right]\\
&\leq&c K E\left[\left(\sup_{s\leq
      t}|x^v_s-x_s|^2\int_0^t|x^v_s-x_s|^2\right)^{1/2}\right]\\
&\leq&\frac{caK}{2}E\left[\sup_{s\leq
      t}|x^v_s-x_s|^2\right]+\frac{cK}{2a}E\int_0^t|x^v_s-x_s|^2ds\,.
\eeaa
Hence we get
\beaa
E\left[\sup_{s\leq t}|x^v_s-x_s|^2\right]&\leq& acKE\left[\sup_{s\leq
    t}|x^v_s-x_s|^2\right]+\frac{cK}{a}E\int_0^t|x^v_s-x_s|^2ds\\
&&+(2K+K^2+\delta)E\int_0^t|x^v_s-x_s|^2ds+\frac{1}{\delta}E\int_0^t|\dot{v}_s|^2ds\,,
\eeaa
where $a,\delta>0$ are arbitrary, $c$ is the constant of Davis'
inequality. From above, we obtain
\beaa
(1-acK)E\left[\sup_{s\leq
    t}|x^v_s-x_s|^2\right]&\leq&\left(\frac{cK}{a}+2K+K^2+\delta\right)E\int_0^t|x^v_s-x_s|^2ds\\
&&+\frac{2}{\delta}H_t(Q|P)
\eeaa
and Gronwall lemma implies that 
\beaa
E\left[\sup_{s\leq
    t}|x^v_s-x_s|^2\right]&\leq&\frac{2}{\delta(1-acK)}H_t(Q|P)\\
&&\cdot
\exp\left[\frac{t}{1-acK}\left(\frac{cK}{a}+\delta+2K+K^2\right)\right]\,.
\eeaa
Taking $t=1$ and minimizing the r.h.s. of the last inequality w.r. to
$\delta$ completes the proof.
\nqed

\noindent
It is important to notice that we did not use any regularity property
about $b$ except that the integrability of $t\to b(t,x)$ for almost
all $x$ in  an intermediate step. This observation means that we can deal with very singular
drifts provided that they are dissipative. Let us give an application of
Proposition  \ref{prop-1} to multi-valued SDE (cf. \cite{CEP}) from this point of view

\begin{theorem}
\label{thm-2}
Let $P$ be the law of the process which is the solution of the
following multi-valued stochastic differential equation:
$$
m(X_t)dt+\sigma (t,X_t)dW_t\in dX_t+A(X_t)dt\,,\, X_0=x\in D(A)\,,
$$
where $A$ is a maximal, monotone set-valued function (hence $-A$ is
dissipative), such that $Int(D(A))\neq \emptyset$. Assume that
$\sigma$ and $m$ are uniformly $K$-Lipschitz and that $\sigma$ is
bounded. Then $P$ satisfies the transportation cost inequality
(\ref{tr-2}). If $\sigma$ is only Lipschitz, but not necessarily
bounded, then $P$ satisfies the inequality (\ref{ineq-3}).

\end{theorem}
\nproof
Let $b_n$ be the Yosida approximation of $A$, i.e.,
$J_n=(I_{\R^d}+\frac{1}{n}A)^{-1}$  and $-b_n=n(I-J_n)$ then $b_n$ is dissipative and
Lipschitz, hence the law of the solution of the SDE
$$
dX^n_t=\sigma (t,X_t)dW_t+b_n(X^n_t)dt+m(X^n_t)dt
$$ 
satisfies the inequality (\ref{tr-2}) with the constants independent
of $n$, moreover the law of $(X^n, n\in
\N)$ converges weakly to $P$ (cf. \cite{CEP}), hence $P$ satisfies
also the inequality (\ref{tr-2}) due to Lemma \ref{stab-lemma}.
\nqed

\noindent
As an example of application of this theorem, let us give
\begin{theorem}
\label{th-3}
Let $P$ be the law of the solution of the following SDE:
$$
dX_t^i=m(X^i_t)dt+\sigma(X_t^i)dW_t^i+\ga\sum_{1\leq j\neq i\leq
  d}\frac{1}{X_t^i-X^j_t}dt\,,\,i=1,\ldots,d\,, 
$$
with $\sigma$ bounded and Lipschitz, $\ga>0$. Then $P$ satisfies the
transportation cost inequality  (\ref{tr-2}) and if $\sigma$ is not
bounded but only Lipschitz, then $P$ satisfies the inequality  (\ref{ineq-3}).
\end{theorem}
\nproof
It suffices to remark that the drift term following $\ga$ is the
subdifferential of the concave function defined by 
$$
F(x)=\ga\sum_{i<j}\log(x^j-x^i)
$$
if $x^1<x^2<\ldots <x^d$ and it is equal to $-\infty$ otherwise.
\nqed

\remark For details about the equation of Theorem \ref{th-3}
cf. \cite{cepa-lep1}. Moreover Theorem \ref{thm-2} is applicable to  all the models given in \cite{cepa-lep2}.

\vspace{2cm}
\footnotesize{
\noindent
A. S. \"Ust\"unel, Institut Telecom, Telecom ParisTech, LTCI CNRS D\'ept. Infres, \\
46, rue Barrault, 75013, Paris, France\\
ustunel@enst.fr}

\end{document}